\documentclass[12pt]{article}
\usepackage{amsmath}
\usepackage{amssymb}
\usepackage{bbm}
\usepackage[margin=2.25cm]{geometry}
\usepackage{hyperref}
\usepackage{tikz}
\usepackage{setspace}

\usetikzlibrary { decorations.pathmorphing, decorations.pathreplacing, decorations.shapes, }

\newtheorem{thm}{Theorem}[section]

\newtheorem{prop}[thm]{Proposition}

\newtheorem{theorem}[thm]{Theorem}
\newtheorem{remark}[thm]{Remark}

\newtheorem{corollary}[thm]{Corollary}

\newtheorem{conj}[thm]{Conjecture}

\newenvironment{proofof}{}{\hfill$\square$\vskip.5cm}

\newcommand{\N}{\mathbb{N}}
\newcommand{\C}{\mathbb{C}}

\newcommand{\Z}{\mathbb{Z}

}

\newcommand{\x}{\boldsymbol{x}}
\newcommand{\y}{\boldsymbol{y}}

\renewcommand{\Z}{{\mathbb{Z}}}

\begin{document} 

\title{Some Observations about the ``Generalized Abundancy Index''}
\author{Shannon Starr}
\date{\today}

\maketitle

\abstract{Let $\mathcal{A}(\ell,n) \subset S_n^{\ell}$ denote the set of all  $\ell$-tuples $(\pi_1,\dots,\pi_{\ell})$, for $\pi_1,\dots,\pi_{\ell} \in S_n$
satisfying: $\forall i<j$ we have $\pi_i\pi_j=\pi_j\pi_i$.
Considering the action of $S_n$ on $[n]=\{1,\dots,n\}$, let 
$\kappa(\pi_1,\dots,\pi_{\ell})$ be equal to the number of  orbits of the action of the subgroup $\langle \pi_1,\dots,\pi_{\ell} \rangle \subset S_n$.
There has been interest in the study of the combinatorial numbers $A(\ell,n,k)$ equal
to the cardinalities $|\{(\pi_1,\dots,\pi_{\ell}) \in \mathcal{A}(\ell,n)\, :\, \kappa(\pi_1,\dots\pi_{\ell})=k\}|$.
If one defines $B(\ell,n)=A(\ell,n,1)/(n-1)!$, then it is known that $B(\ell,n) = \sum_{(f_1,\dots,f_{\ell}) \in \N^{\ell}} \mathbf{1}_{\{n\}}(f_1\cdots f_{\ell}) \prod_{r=1}^{\ell-1} f_r^{\ell-r}$.
A special case, $\ell=2$, is $B(2,n) = \sum_{d|n} d = \sigma_1(n)$ the sum-of-divisors function.
Then $A(2,n,1)/n!=B(2,n)/n$ is called the abundancy index: $\sigma_1(n)/n$.
We call $B(\ell,n) n^{-\ell+1}$ the ``generalized abundancy index.''
Building on work of Abdesselam, using the probability model, we prove that $\lim_{N \to \infty} N^{-1} \sum_{n=1}^{N} B(\ell,n) n^{-\ell+1}$ equals
$\zeta(2)\cdots \zeta(\ell)$.
Motivated by this we state a more precise conjecture
for the asymptotics of $-\zeta(2) + N^{-1}\sum_{n=1}^{N} (B(2,n)/n)$.
}

\section{Introduction}

Let $S_n$ denote the group of permutations of $[n]=\{1,\dots,n\}$, viewed as bijections of $[n]$ to itself, so that 
we write $\pi(j)$  for the image of the point $j$ under the action of $\pi \in S_n$.
Let $S_n^{\ell}$ denote the set of all $\ell$-tuples $(\pi_1,\dots,\pi_{\ell})$, with $\pi_1,\dots,\pi_{\ell} \in S_n$.
Finally, let $\mathcal{A}(\ell,n)$ denote the set of all $\ell$-tuples $(\pi_1,\dots,\pi_{\ell}) \in S_n^{\ell}$
satisfying
$$
\forall 1\leq i<j\leq \ell\, ,\ \text{ we have }\ \pi_i \circ \pi_{j} = \pi_j \circ \pi_i\, .
$$
There is interest in the following combinatorial numbers
$$
	A(\ell,n,k)\, =\, |\{(\pi_1,\dots,\pi_{\ell}) \in \mathcal{A}(\ell,n)\, :\, \kappa(\pi_1,\dots,\pi_{\ell}) = k\}|\, .
$$
Letting $\langle \pi_1,\dots,\pi_{\ell} \rangle$ denote the subgroup of $S_n$ generated by $\pi_1,\dots,\pi_{\ell}$, we let $\kappa(\pi_1,\dots,\pi_{\ell})$ denote the number of orbits
of $\langle \pi_1,\dots,\pi_{\ell} \rangle$ as it acts on $[n]$: defining
$$
\langle \pi_1,\dots,\pi_{\ell} \rangle\{j\}\, =\, \{\rho(j)\, :\, \rho \in \langle \pi_1,\dots,\pi_{\ell} \rangle\}\, ,
$$
for each $j \in [n]$,
then the number of orbits is 
$$
\kappa(\pi_1,\dots,\pi_{\ell})\, =\, |\{\langle \pi_1,\dots,\pi_{\ell} \rangle\{j\}\, :\, j \in [n]\}|\, .
$$
Let $\mathcal{A}(\ell,n,k) = \{(\pi_1,\dots,\pi_{\ell}) \in \mathcal{A}(\ell,n)\, :\, 
\kappa(\pi_1,\dots,\pi_{\ell})=k\}$.
It is not hard to see that
\begin{equation}
\label{eq:Adecomp}
A(\ell,n,k)\, =\, \frac{n!}{k!} \sum_{\substack{(\nu_1,\dots,\nu_k) \in \N^k \\ \nu_1+\dots+\nu_k=n}}\, \prod_{r=1}^{k} \left(\frac{A(\ell,\nu_k,1)}{\nu_k!}\right)\, .
\end{equation}
We note that the notion of the Bell transform is an important concept occurring
frequently in combinatorics and probability.
We suggest Chapter 1 of Pitman's excellent textbook \cite{Pitman}.
Motivated by this, one may wish to study the combinatorial numbers $A(\ell,n,1)$.
In fact, the more common definition is 
$$
	B(\ell,n)\, =\, \frac{A(\ell,n,1)}{(n-1)!}\, .
$$
One may easily see, by induction, that
\begin{equation}
\label{eq:Abdesselam}
	B(\ell,n)\, =\, \sum_{\substack{(f_1,\dots,f_{\ell}) \in \N^{\ell}\\ f_1\cdots f_{\ell}=n}}
	f_1^{\ell-1} f_2^{\ell-2} \cdots f_{\ell-1}\, ,
\end{equation}
which is a formula noted in a recent article of Abdesselam \cite{Abdesselam}.
Abdesselam notes a different formula, too.
Let $d_1|d_2|\dots|d_{\ell-1}|n$ denote a multiplicative flag $(d_1,\dots,d_{\ell-1}) \in \N^{\ell-1}$ satisfying
$$
d_1|d_2\, ,\ d_2|d_3\, ,\ \dots\, ,\ d_{\ell-1}|n\, .
$$
Then
$$
B(\ell,n)\, =\, \sum_{d_1|d_2|\dots|d_{\ell-1}|n} d_1d_2\cdots d_{\ell-1}\, .
$$
Clearly, we could write $d_r = f_1\cdots f_r$ for $r=1,\dots,\ell-1$ to map between these two.
A well-known result of Bryan and Fulman \cite{BryanFulman} is
$$
	1+\sum_{n=1}^{\infty} \sum_{k=0}^{n} \frac{A(\ell,n,k)}{n!} x^k z^n\,
	=\, \exp\left(-x\sum_{(f_1,\dots,f_{\ell-1}) \in \N^{\ell-1}}f_1^{\ell-2} f_2^{\ell-3} \cdots f_{\ell-2}
\ln(1-z^{f_1\cdots f_{\ell}})\right)\, .
$$
In particular, that means
\begin{equation}
\begin{split}
	\sum_{n=1}^{\infty} \frac{B(\ell,n)}{n}\, z^n\, 
	&=\, -\sum_{(f_1,\dots,f_{\ell-1}) \in \N^{\ell-1}}f_1^{\ell-2} f_2^{\ell-3} \cdots f_{\ell-2}
\ln(1-z^{f_1\cdots f_{\ell-1}})\\
	&=\, \sum_{(f_1,\dots,f_{\ell}) \in \N^{\ell}} z^{f_1\cdots f_{\ell}} \prod_{r=1}^{\ell} f_r^{\ell-1-r}\, ,
\end{split}
\end{equation}
which implies (\ref{eq:Abdesselam}).
In \cite{AbdesselamStudents}, they gave a re-derivation of the equation (\ref{eq:Abdesselam}).
We will also give a description of this in 
Appendix \ref{app:tori}.

We note that for many purposes, the fractions $A(\ell,n,k)/n!$ are fundamental objects.
For instance, the Cauchy integral formula gives
\begin{equation}
\label{eq:Cauchy}
	\frac{A(\ell,n,k)}{n!}\, =\, e^{-\ln(k!)} \oint_{C(0;r)} e^{-n\ln(r)-in\theta+k \ln\left(\sum_{\nu\in\N} r^n e^{in\theta} (B(\ell,\nu)/\nu)\right)}\, \frac{d\theta}{2\pi}\, ,
\end{equation}
for any $r \in (0,1)$.

Abdesselam also noted that the numbers $B(\ell,n)$ are multiplicative, meaning that if 
$$
n\, =\, p_1^{a(1)} p_2^{a(2)}\cdots\, ,
$$
where $p_1,p_2,\dots$ is the ordered enumeration of primes, and $a(1),a(2),\dots \in \{0,1,\dots\}$,
then 
$$
B(\ell,n)\, =\, \prod_{r=1}^{\infty} \widehat{B}(\ell,p_r,a(r))\, ,
$$
for a function which Abdesselam also calculated
$$
\widehat{B}(\ell,p,a)\, =\, p^{(\ell-1)a}\, \frac{(p^{-a-1};p^{-1})_{\ell-1}}{(p^{-1};p^{-1})_{\ell-1}}\, ,
$$
where the $q$-Pocchammer symbol is 
$$
	(z;q)_r\, =\, (1-z)(1-qz)(1-q^2z)\cdots(1-q^{r-1}z)\, ,
$$
for $r \in \{0,1,2,\dots\}$ as usual.
Using this, it is easy to prove the following.
\begin{theorem}
\label{thm:Cesaro}
For any $\ell \in \{2,3,\dots\}$, and
$|q|<1$, we have a $q$-integral analogue of the power rule:
$$
	z(1-q)\sum_{k=0}^{\infty}q^{k} (q^{k+1}z;q)_{\ell-1} \, 
	=\, \frac{(1-q)\left(1- (z;q)_{\ell}\right)}{1-q^{\ell}}\, .
$$
Specializing this to $z=1$, and using the probability model, we have
$$
	\lim_{N \to \infty} \frac{1}{N}\, \sum_{n=1}^{N} \frac{B(\ell,n)}{n^{\ell-1}}\, 
	=\, \zeta(2)\cdots \zeta(\ell)\, .
$$
\end{theorem}
We will prove this in Appendix \ref{sec:q}.
The first part of the result is just an application of the $q$-binomial formula.
The idea for the second part is to use the probability model, as given for example in the excellent textbook 
of Arratia, Tavar\'e and Barbour \cite{ATB}.
We note that this is Cesaro sum convergence of the fractions $B(\ell,n)/n^{\ell-1}$.
That is because $B(\ell,n)/n^{\ell-1}$ does not converge directly for $\ell \in \{2,3,\dots\}$.
Instead, it has a limiting distribution.
The limiting distribution may also be understood
using the probability model and Abdesselam's formula.
But the formulas are more involved. For example, we may view $\zeta(2)\cdots \zeta(\ell)$ as the first moment of the limiting distribution of $B(\ell,n)/n^{\ell-1}$.
Then one may calculate the analogue of the second moment in the case that $\ell=2$:
$\zeta(2)^2\zeta(3)/\zeta(4)$.
This is as explicit as the first moment.
But the second moment in the case $\ell=3$ is
$$
	\zeta(2)\zeta(3)\zeta(4)\zeta(5) \prod_{n=1}^{\infty} \left(1+\frac{1}{p_n}+\frac{2}{p_n^3}+\frac{1}{p_n^4}+\frac{1}{p_n^6}\right)\, .
$$
The third moment in the case $\ell=2$ is  $\zeta(3)\zeta(4) \prod_{n=1}^{\infty} \left(1-\frac{1}{p_n^3}+\frac{3}{p_n(p_n-1)}\right)$.
These formulas are not explicit.
\begin{remark}
To better understand the particular choice of $B(\ell,n)/n^{\ell-1}$, use
the convention of Abdesselam:
$$
\frac{B(\ell,n)}{n^{\ell-1}}\, =\, \sum_{d_1|\dots|d_{\ell-1}|n}\, \prod_{r=1}^{\ell-1}
\left(\frac{d_r}{n}\right)\, .
$$
\end{remark}

\section{Review of rigorous results}

Suppose that one writes
\begin{equation}
\label{eq:HcalDef}
\mathcal{H}_{\ell,n}(x)\, =\, \sum_{k=0}^{n} \frac{1}{n!}\, A(\ell,n,k) x^k\, ,
\end{equation}
for $x \in \C$.
This may already be viewed as a generating function.
Define the two-variable generating function
$$
\mathcal{G}_{\ell}(x,z)\, =\, \sum_{n=0}^{\infty} \mathcal{H}_{\ell,n}(x) z^n\, .
$$
Then Bryan and Fulman showed
\begin{equation}
\label{eq:BryanFulman}
	\mathcal{G}_{\ell}(x,z)\, =\, \exp\left(x \mathcal{L}_{\ell}(z)\right)\, ,
\end{equation}
for $|z|<1$, where
$$
	\mathcal{L}_{\ell}(z)\, =\, -\sum_{\delta_1,\dots,\delta_{\ell-1} \in \N} \delta_1^{\ell-2} \delta_2^{\ell-3} \cdots \delta_{\ell-2}\, \ln\left(1-z^{\delta_1 \cdots \delta_{\ell-1}}\right)\, .
$$
Using this and the Hardy-Ramanujan circle method, and the Euler-Maclaurin sum-to-integral comparison formula, one has the following.
\begin{prop}
For each $x \in (0,\infty)$ and $\ell \in \{2,3,\dots\}$, we have
$$
\ln\left(\mathcal{H}_{\ell,n}(x)\right)\, \sim\, n^{(\ell-1)/\ell}\, \cdot 
\frac{\ell}{\ell-1}\, \Big(x\mathcal{K}_{\ell}\Big)^{1/\ell}\, ,\ \text{ as $n \to \infty$,}
$$
where
$$
\mathcal{K}_{\ell}\, =\, \Gamma(\ell) \zeta(2)\cdots \zeta(\ell)\, .
$$
\end{prop}
This is the subject of a forthcoming preprint of Abdesselam and the author. The Legendre transform leads to a large deviation theorem for the probability mass function $k \in \{0,\dots,n\} \mapsto A(\ell,n,k) x^k /n!$, in the large $n$ limit.
A corollary of this is a central limit theorem, which we prove in our forthcoming preprint.
This may be stated as follows:
\begin{corollary}
\label{cor:CLT}
For any $y \in (0,\infty)$ and $-\infty<\alpha<\beta<\infty$, we have,
for $\mathcal{I}_n(y;\alpha,\beta) = \{\lfloor y n^{(\ell-1)/\ell}+\alpha n^{(\ell-1)/(2\ell)}\rfloor\, ,\ \dots\, ,\ \lfloor y n^{(\ell-1)/\ell}+\beta n^{(\ell-1)/(2\ell)}\rfloor\}$,
$$
\sum_{k\in \mathcal{I}_n(y;\alpha,\beta)}
	\frac{A(\ell,n,k)}{n!}\, 
=\, 
\sum_{k\in \mathcal{I}_n(y;\alpha,\beta)}
\frac{\sqrt{\ell-1}}{2\pi n}\, 
e^{k\ell - k\ell \ln(k) + (\ell-1)k\ln(n)-\lambda_{\ell} k} \left(1+o(1)\right)\, ,
$$
where $\lambda_{\ell} = \ell \ln(\ell-1) - \ln(\mathcal{K}_{\ell})$.
\end{corollary}

However, one may make the following conjecture, which appears to be true.
\begin{conj}
\label{conj:LCLT}
For $n \to \infty$ and $k \to \infty$ such that $k/n^{(\ell-1)/\ell}$ converges to a quantity in $(0,\infty)$, we should have
$$
\frac{A(\ell,n,k)}{n!}\, =\, \frac{\sqrt{\ell-1}}{2\pi n}\, e^{k\ell-k\ell\, \ln(k) + (\ell-1) k\ln(n) - \widetilde{\lambda}_{\ell}\left(\frac{k}{n}\right) k}(1+o(1))\, ,
$$
where
$$
\widetilde{\lambda}_{\ell}\left(\frac{k}{n}\right)\, =\, \lambda_{\ell} + \varepsilon_{\ell}\left(\frac{k}{n}\right)\, ,
$$
where $\varepsilon_{\ell}(t)$ is a linear combination of terms of the form $t^j \left(\ln(t)\right)^k$ for $1\leq j\leq \ell-1$ (and $k \in \{0,1,\dots\}$).
\end{conj}
The reasons for the corrections to $\lambda_{\ell}$ represented in $\widetilde{\lambda}_{\ell}(k/n)$ is that the CLT stops at a certain order
which is not small enough to get to order-1.
Moreover, for $\mathcal{H}_{\ell,n}(1)$, a recent paper of Bringmann, Franke and Heim \cite{BFH} has indicated how to calculate corrections
down to order 1.

As an example, for $\ell=2$, we have
$$
\mathcal{H}_{2,n}(1)\, =\, p(n)\, ,
$$
the partition number. Then Hardy and Ramanujan's formula gives
\begin{equation}
\label{eq:varepsilon}
\varepsilon_2(t)\, =\, \frac{\ln(24 \zeta(2))}{4\zeta(2)}\, t - \frac{1}{2\zeta(2)}\, t \ln(t)\, ,
\end{equation}
so that $\mathcal{H}_{2,n}(x) = x^{(1+x)/4} 2^{-(5+3x)/4} 3^{-(1+x)/4} n^{-(3+x)/4} \exp(2\sqrt{x\zeta(2)n})$
times $1+o(1)$ as $n \to \infty$.
Using this we may make the more refined conjecture for the abundancy index.
\begin{conj}[Main conjecture]
\label{conj:main}
We have
$$
\lim_{\mathcal{N} \to \infty} 
\frac{1}{\mathcal{N}} \sum_{N=1}^{\mathcal{N}} 
\left(\sum_{n=1}^{N} \frac{B(2,n)}{n}-\zeta(2) N+\frac{1}{2}\, \ln(N)\right)\, =\, -\mu\, ,
$$
for a constant $\mu$ expressible in terms of $\gamma$, the Euler-Mascheroni constant: 
$$
\mu\, =\, \frac{\gamma}{2} + \frac{\ln(24 \zeta(2))}{4} - \frac{\zeta(2)}{2}\, \approx\, 0.38507933223132607\, .
$$
\end{conj}
In section 4, we will explain why guessing Conjecture \ref{conj:LCLT}
leads to guessing Conjecture \ref{conj:main}.
First we present numerical evidence.

\section{Python algorithm and numerical evidence}

A Python algorithm which generates the abundancy index for numbers
up to $N=10^6$ in a time about 1 hour is this:
\begin{verbatim}
import numpy as np
import sympy
import pickle

N_max = 1000000
B1_lst = np.array([])
for n_ctr in range(N_max):
    n_num = n_ctr+1
    d_lst = sympy.divisors(n_num)
    B1_num = np.sum(d_lst)
    B1_lst = np.append(B1_lst,B1_num)

abndcy_lst = np.divide(B1_lst,1+np.array(range(N_max)))

with open('abndcy.pickle', 'wb') as f:
    # Pickle the 'data' dictionary using the highest protocol available.
    pickle.dump(abndcy_lst, f, pickle.HIGHEST_PROTOCOL)
\end{verbatim}

Then another short Python code to determine the sample mean of the error data is this:
\begin{verbatim}
import numpy as np
import pickle

with open('abndcy.pickle', 'rb') as f:
    # The protocol version used is detected automatically, so we do not
    # have to specify it.
    abndcy_lst = pickle.load(f)

N_max=1000000
first_approx = np.pi**2/6*(1+np.array(range(N_max)))
second_approx = -0.5*np.log(1+np.array(range(N_max)))
err_data = np.cumsum(abndcy_lst)-first_approx-second_approx
np.mean(err_data)
\end{verbatim}
Then output is 
$$
-0.38508487292161986
$$
So the relative error is approximately $0.0014\%$.

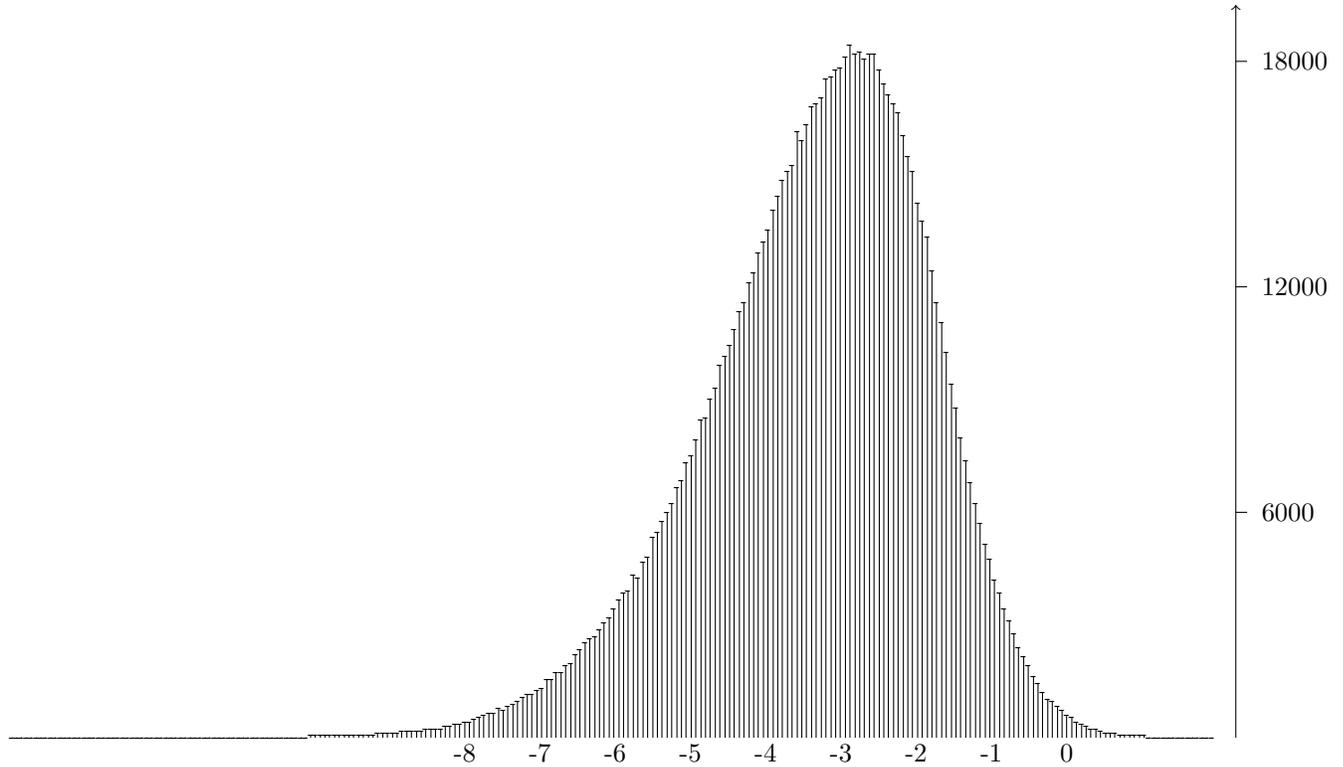
\begin{figure}
{\tiny
$$
\begin{tikzpicture}[xscale=1,yscale=3]
\draw ( -14.02 , 0 )--( -14.02 , 0.0 ) node[] {-};
\draw ( -13.96 , 0 )--( -13.96 , 0.0 ) node[] {-};
\draw ( -13.9 , 0 )--( -13.9 , 0.0 ) node[] {-};
\draw ( -13.83 , 0 )--( -13.83 , 0.0 ) node[] {-};
\draw ( -13.77 , 0 )--( -13.77 , 0.0 ) node[] {-};
\draw ( -13.7 , 0 )--( -13.7 , 0.0 ) node[] {-};
\draw ( -13.64 , 0 )--( -13.64 , 0.0 ) node[] {-};
\draw ( -13.58 , 0 )--( -13.58 , 0.0 ) node[] {-};
\draw ( -13.51 , 0 )--( -13.51 , 0.0 ) node[] {-};
\draw ( -13.45 , 0 )--( -13.45 , 0.0 ) node[] {-};
\draw ( -13.38 , 0 )--( -13.38 , 0.0 ) node[] {-};
\draw ( -13.32 , 0 )--( -13.32 , 0.0 ) node[] {-};
\draw ( -13.26 , 0 )--( -13.26 , 0.0 ) node[] {-};
\draw ( -13.19 , 0 )--( -13.19 , 0.0 ) node[] {-};
\draw ( -13.13 , 0 )--( -13.13 , 0.0 ) node[] {-};
\draw ( -13.06 , 0 )--( -13.06 , 0.0 ) node[] {-};
\draw ( -13.0 , 0 )--( -13.0 , 0.0 ) node[] {-};
\draw ( -12.94 , 0 )--( -12.94 , 0.0 ) node[] {-};
\draw ( -12.87 , 0 )--( -12.87 , 0.0 ) node[] {-};
\draw ( -12.81 , 0 )--( -12.81 , 0.0 ) node[] {-};
\draw ( -12.74 , 0 )--( -12.74 , 0.0 ) node[] {-};
\draw ( -12.68 , 0 )--( -12.68 , 0.0 ) node[] {-};
\draw ( -12.62 , 0 )--( -12.62 , 0.0 ) node[] {-};
\draw ( -12.55 , 0 )--( -12.55 , 0.0 ) node[] {-};
\draw ( -12.49 , 0 )--( -12.49 , 0.0 ) node[] {-};
\draw ( -12.42 , 0 )--( -12.42 , 0.0 ) node[] {-};
\draw ( -12.36 , 0 )--( -12.36 , 0.0 ) node[] {-};
\draw ( -12.3 , 0 )--( -12.3 , 0.0 ) node[] {-};
\draw ( -12.23 , 0 )--( -12.23 , 0.0 ) node[] {-};
\draw ( -12.17 , 0 )--( -12.17 , 0.0 ) node[] {-};
\draw ( -12.1 , 0 )--( -12.1 , 0.0 ) node[] {-};
\draw ( -12.04 , 0 )--( -12.04 , 0.0 ) node[] {-};
\draw ( -11.97 , 0 )--( -11.97 , 0.0 ) node[] {-};
\draw ( -11.91 , 0 )--( -11.91 , 0.0 ) node[] {-};
\draw ( -11.85 , 0 )--( -11.85 , 0.0 ) node[] {-};
\draw ( -11.78 , 0 )--( -11.78 , 0.0 ) node[] {-};
\draw ( -11.72 , 0 )--( -11.72 , 0.0 ) node[] {-};
\draw ( -11.65 , 0 )--( -11.65 , 0.0 ) node[] {-};
\draw ( -11.59 , 0 )--( -11.59 , 0.0 ) node[] {-};
\draw ( -11.53 , 0 )--( -11.53 , 0.0 ) node[] {-};
\draw ( -11.46 , 0 )--( -11.46 , 0.0 ) node[] {-};
\draw ( -11.4 , 0 )--( -11.4 , 0.0 ) node[] {-};
\draw ( -11.33 , 0 )--( -11.33 , 0.0 ) node[] {-};
\draw ( -11.27 , 0 )--( -11.27 , 0.0 ) node[] {-};
\draw ( -11.21 , 0 )--( -11.21 , 0.0 ) node[] {-};
\draw ( -11.14 , 0 )--( -11.14 , 0.0 ) node[] {-};
\draw ( -11.08 , 0 )--( -11.08 , 0.0 ) node[] {-};
\draw ( -11.01 , 0 )--( -11.01 , 0.0 ) node[] {-};
\draw ( -10.95 , 0 )--( -10.95 , 0.0 ) node[] {-};
\draw ( -10.89 , 0 )--( -10.89 , 0.0 ) node[] {-};
\draw ( -10.82 , 0 )--( -10.82 , 0.0 ) node[] {-};
\draw ( -10.76 , 0 )--( -10.76 , 0.0 ) node[] {-};
\draw ( -10.69 , 0 )--( -10.69 , 0.0 ) node[] {-};
\draw ( -10.63 , 0 )--( -10.63 , 0.0 ) node[] {-};
\draw ( -10.57 , 0 )--( -10.57 , 0.0 ) node[] {-};
\draw ( -10.5 , 0 )--( -10.5 , 0.0 ) node[] {-};
\draw ( -10.44 , 0 )--( -10.44 , 0.0 ) node[] {-};
\draw ( -10.37 , 0 )--( -10.37 , 0.0 ) node[] {-};
\draw ( -10.31 , 0 )--( -10.31 , 0.0 ) node[] {-};
\draw ( -10.25 , 0 )--( -10.25 , 0.0 ) node[] {-};
\draw ( -10.18 , 0 )--( -10.18 , 0.0 ) node[] {-};
\draw ( -10.12 , 0 )--( -10.12 , 0.0 ) node[] {-};
\draw ( -10.05 , 0 )--( -10.05 , 0.01 ) node[] {-};
\draw ( -9.99 , 0 )--( -9.99 , 0.01 ) node[] {-};
\draw ( -9.92 , 0 )--( -9.92 , 0.01 ) node[] {-};
\draw ( -9.86 , 0 )--( -9.86 , 0.01 ) node[] {-};
\draw ( -9.8 , 0 )--( -9.8 , 0.01 ) node[] {-};
\draw ( -9.73 , 0 )--( -9.73 , 0.01 ) node[] {-};
\draw ( -9.67 , 0 )--( -9.67 , 0.01 ) node[] {-};
\draw ( -9.6 , 0 )--( -9.6 , 0.01 ) node[] {-};
\draw ( -9.54 , 0 )--( -9.54 , 0.01 ) node[] {-};
\draw ( -9.48 , 0 )--( -9.48 , 0.01 ) node[] {-};
\draw ( -9.41 , 0 )--( -9.41 , 0.01 ) node[] {-};
\draw ( -9.35 , 0 )--( -9.35 , 0.01 ) node[] {-};
\draw ( -9.28 , 0 )--( -9.28 , 0.01 ) node[] {-};
\draw ( -9.22 , 0 )--( -9.22 , 0.01 ) node[] {-};
\draw ( -9.16 , 0 )--( -9.16 , 0.02 ) node[] {-};
\draw ( -9.09 , 0 )--( -9.09 , 0.02 ) node[] {-};
\draw ( -9.03 , 0 )--( -9.03 , 0.02 ) node[] {-};
\draw ( -8.96 , 0 )--( -8.96 , 0.02 ) node[] {-};
\draw ( -8.9 , 0 )--( -8.9 , 0.02 ) node[] {-};
\draw ( -8.84 , 0 )--( -8.84 , 0.03 ) node[] {-};
\draw ( -8.77 , 0 )--( -8.77 , 0.03 ) node[] {-};
\draw ( -8.71 , 0 )--( -8.71 , 0.03 ) node[] {-};
\draw ( -8.64 , 0 )--( -8.64 , 0.03 ) node[] {-};
\draw ( -8.58 , 0 )--( -8.58 , 0.03 ) node[] {-};
\draw ( -8.52 , 0 )--( -8.52 , 0.04 ) node[] {-};
\draw ( -8.45 , 0 )--( -8.45 , 0.04 ) node[] {-};
\draw ( -8.39 , 0 )--( -8.39 , 0.04 ) node[] {-};
\draw ( -8.32 , 0 )--( -8.32 , 0.04 ) node[] {-};
\draw ( -8.26 , 0 )--( -8.26 , 0.05 ) node[] {-};
\draw ( -8.2 , 0 )--( -8.2 , 0.05 ) node[] {-};
\draw ( -8.13 , 0 )--( -8.13 , 0.06 ) node[] {-};
\draw ( -8.07 , 0 )--( -8.07 , 0.06 ) node[] {-};
\draw ( -8.0 , 0 )--( -8.0 , 0.07 ) node[] {-};
\draw ( -7.94 , 0 )--( -7.94 , 0.07 ) node[] {-};
\draw ( -7.88 , 0 )--( -7.88 , 0.08 ) node[] {-};
\draw ( -7.81 , 0 )--( -7.81 , 0.09 ) node[] {-};
\draw ( -7.75 , 0 )--( -7.75 , 0.1 ) node[] {-};
\draw ( -7.68 , 0 )--( -7.68 , 0.11 ) node[] {-};
\draw ( -7.62 , 0 )--( -7.62 , 0.11 ) node[] {-};
\draw ( -7.55 , 0 )--( -7.55 , 0.13 ) node[] {-};
\draw ( -7.49 , 0 )--( -7.49 , 0.12 ) node[] {-};
\draw ( -7.43 , 0 )--( -7.43 , 0.14 ) node[] {-};
\draw ( -7.36 , 0 )--( -7.36 , 0.15 ) node[] {-};
\draw ( -7.3 , 0 )--( -7.3 , 0.16 ) node[] {-};
\draw ( -7.23 , 0 )--( -7.23 , 0.18 ) node[] {-};
\draw ( -7.17 , 0 )--( -7.17 , 0.19 ) node[] {-};
\draw ( -7.11 , 0 )--( -7.11 , 0.19 ) node[] {-};
\draw ( -7.04 , 0 )--( -7.04 , 0.21 ) node[] {-};
\draw ( -6.98 , 0 )--( -6.98 , 0.22 ) node[] {-};
\draw ( -6.91 , 0 )--( -6.91 , 0.26 ) node[] {-};
\draw ( -6.85 , 0 )--( -6.85 , 0.26 ) node[] {-};
\draw ( -6.79 , 0 )--( -6.79 , 0.29 ) node[] {-};
\draw ( -6.72 , 0 )--( -6.72 , 0.29 ) node[] {-};
\draw ( -6.66 , 0 )--( -6.66 , 0.32 ) node[] {-};
\draw ( -6.59 , 0 )--( -6.59 , 0.33 ) node[] {-};
\draw ( -6.53 , 0 )--( -6.53 , 0.37 ) node[] {-};
\draw ( -6.47 , 0 )--( -6.47 , 0.39 ) node[] {-};
\draw ( -6.4 , 0 )--( -6.4 , 0.42 ) node[] {-};
\draw ( -6.34 , 0 )--( -6.34 , 0.44 ) node[] {-};
\draw ( -6.27 , 0 )--( -6.27 , 0.45 ) node[] {-};
\draw ( -6.21 , 0 )--( -6.21 , 0.48 ) node[] {-};
\draw ( -6.15 , 0 )--( -6.15 , 0.51 ) node[] {-};
\draw ( -6.08 , 0 )--( -6.08 , 0.53 ) node[] {-};
\draw ( -6.02 , 0 )--( -6.02 , 0.57 ) node[] {-};
\draw ( -5.95 , 0 )--( -5.95 , 0.61 ) node[] {-};
\draw ( -5.89 , 0 )--( -5.89 , 0.64 ) node[] {-};
\draw ( -5.83 , 0 )--( -5.83 , 0.65 ) node[] {-};
\draw ( -5.76 , 0 )--( -5.76 , 0.72 ) node[] {-};
\draw ( -5.7 , 0 )--( -5.7 , 0.71 ) node[] {-};
\draw ( -5.63 , 0 )--( -5.63 , 0.78 ) node[] {-};
\draw ( -5.57 , 0 )--( -5.57 , 0.8 ) node[] {-};
\draw ( -5.5 , 0 )--( -5.5 , 0.89 ) node[] {-};
\draw ( -5.44 , 0 )--( -5.44 , 0.91 ) node[] {-};
\draw ( -5.38 , 0 )--( -5.38 , 0.96 ) node[] {-};
\draw ( -5.31 , 0 )--( -5.31 , 1.0 ) node[] {-};
\draw ( -5.25 , 0 )--( -5.25 , 1.04 ) node[] {-};
\draw ( -5.18 , 0 )--( -5.18 , 1.11 ) node[] {-};
\draw ( -5.12 , 0 )--( -5.12 , 1.14 ) node[] {-};
\draw ( -5.06 , 0 )--( -5.06 , 1.22 ) node[] {-};
\draw ( -4.99 , 0 )--( -4.99 , 1.25 ) node[] {-};
\draw ( -4.93 , 0 )--( -4.93 , 1.32 ) node[] {-};
\draw ( -4.86 , 0 )--( -4.86 , 1.41 ) node[] {-};
\draw ( -4.8 , 0 )--( -4.8 , 1.42 ) node[] {-};
\draw ( -4.74 , 0 )--( -4.74 , 1.5 ) node[] {-};
\draw ( -4.67 , 0 )--( -4.67 , 1.55 ) node[] {-};
\draw ( -4.61 , 0 )--( -4.61 , 1.65 ) node[] {-};
\draw ( -4.54 , 0 )--( -4.54 , 1.69 ) node[] {-};
\draw ( -4.48 , 0 )--( -4.48 , 1.74 ) node[] {-};
\draw ( -4.42 , 0 )--( -4.42 , 1.81 ) node[] {-};
\draw ( -4.35 , 0 )--( -4.35 , 1.89 ) node[] {-};
\draw ( -4.29 , 0 )--( -4.29 , 1.93 ) node[] {-};
\draw ( -4.22 , 0 )--( -4.22 , 2.02 ) node[] {-};
\draw ( -4.16 , 0 )--( -4.16 , 2.06 ) node[] {-};
\draw ( -4.1 , 0 )--( -4.1 , 2.15 ) node[] {-};
\draw ( -4.03 , 0 )--( -4.03 , 2.2 ) node[] {-};
\draw ( -3.97 , 0 )--( -3.97 , 2.25 ) node[] {-};
\draw ( -3.9 , 0 )--( -3.9 , 2.34 ) node[] {-};
\draw ( -3.84 , 0 )--( -3.84 , 2.4 ) node[] {-};
\draw ( -3.78 , 0 )--( -3.78 , 2.47 ) node[] {-};
\draw ( -3.71 , 0 )--( -3.71 , 2.51 ) node[] {-};
\draw ( -3.65 , 0 )--( -3.65 , 2.54 ) node[] {-};
\draw ( -3.58 , 0 )--( -3.58 , 2.69 ) node[] {-};
\draw ( -3.52 , 0 )--( -3.52 , 2.65 ) node[] {-};
\draw ( -3.46 , 0 )--( -3.46 , 2.72 ) node[] {-};
\draw ( -3.39 , 0 )--( -3.39 , 2.8 ) node[] {-};
\draw ( -3.33 , 0 )--( -3.33 , 2.81 ) node[] {-};
\draw ( -3.26 , 0 )--( -3.26 , 2.84 ) node[] {-};
\draw ( -3.2 , 0 )--( -3.2 , 2.92 ) node[] {-};
\draw ( -3.13 , 0 )--( -3.13 , 2.93 ) node[] {-};
\draw ( -3.07 , 0 )--( -3.07 , 2.96 ) node[] {-};
\draw ( -3.01 , 0 )--( -3.01 , 2.97 ) node[] {-};
\draw ( -2.94 , 0 )--( -2.94 , 3.02 ) node[] {-};
\draw ( -2.88 , 0 )--( -2.88 , 3.07 ) node[] {-};
\draw ( -2.81 , 0 )--( -2.81 , 3.03 ) node[] {-};
\draw ( -2.75 , 0 )--( -2.75 , 3.04 ) node[] {-};
\draw ( -2.69 , 0 )--( -2.69 , 3.01 ) node[] {-};
\draw ( -2.62 , 0 )--( -2.62 , 3.03 ) node[] {-};
\draw ( -2.56 , 0 )--( -2.56 , 3.03 ) node[] {-};
\draw ( -2.49 , 0 )--( -2.49 , 2.96 ) node[] {-};
\draw ( -2.43 , 0 )--( -2.43 , 2.9 ) node[] {-};
\draw ( -2.37 , 0 )--( -2.37 , 2.85 ) node[] {-};
\draw ( -2.3 , 0 )--( -2.3 , 2.81 ) node[] {-};
\draw ( -2.24 , 0 )--( -2.24 , 2.77 ) node[] {-};
\draw ( -2.17 , 0 )--( -2.17 , 2.67 ) node[] {-};
\draw ( -2.11 , 0 )--( -2.11 , 2.58 ) node[] {-};
\draw ( -2.05 , 0 )--( -2.05 , 2.51 ) node[] {-};
\draw ( -1.98 , 0 )--( -1.98 , 2.37 ) node[] {-};
\draw ( -1.92 , 0 )--( -1.92 , 2.29 ) node[] {-};
\draw ( -1.85 , 0 )--( -1.85 , 2.22 ) node[] {-};
\draw ( -1.79 , 0 )--( -1.79 , 2.07 ) node[] {-};
\draw ( -1.73 , 0 )--( -1.73 , 1.93 ) node[] {-};
\draw ( -1.66 , 0 )--( -1.66 , 1.84 ) node[] {-};
\draw ( -1.6 , 0 )--( -1.6 , 1.71 ) node[] {-};
\draw ( -1.53 , 0 )--( -1.53 , 1.57 ) node[] {-};
\draw ( -1.47 , 0 )--( -1.47 , 1.46 ) node[] {-};
\draw ( -1.41 , 0 )--( -1.41 , 1.33 ) node[] {-};
\draw ( -1.34 , 0 )--( -1.34 , 1.23 ) node[] {-};
\draw ( -1.28 , 0 )--( -1.28 , 1.13 ) node[] {-};
\draw ( -1.21 , 0 )--( -1.21 , 1.04 ) node[] {-};
\draw ( -1.15 , 0 )--( -1.15 , 0.95 ) node[] {-};
\draw ( -1.08 , 0 )--( -1.08 , 0.86 ) node[] {-};
\draw ( -1.02 , 0 )--( -1.02 , 0.79 ) node[] {-};
\draw ( -0.96 , 0 )--( -0.96 , 0.7 ) node[] {-};
\draw ( -0.89 , 0 )--( -0.89 , 0.64 ) node[] {-};
\draw ( -0.83 , 0 )--( -0.83 , 0.57 ) node[] {-};
\draw ( -0.76 , 0 )--( -0.76 , 0.52 ) node[] {-};
\draw ( -0.7 , 0 )--( -0.7 , 0.46 ) node[] {-};
\draw ( -0.64 , 0 )--( -0.64 , 0.4 ) node[] {-};
\draw ( -0.57 , 0 )--( -0.57 , 0.36 ) node[] {-};
\draw ( -0.51 , 0 )--( -0.51 , 0.32 ) node[] {-};
\draw ( -0.44 , 0 )--( -0.44 , 0.27 ) node[] {-};
\draw ( -0.38 , 0 )--( -0.38 , 0.24 ) node[] {-};
\draw ( -0.32 , 0 )--( -0.32 , 0.2 ) node[] {-};
\draw ( -0.25 , 0 )--( -0.25 , 0.17 ) node[] {-};
\draw ( -0.19 , 0 )--( -0.19 , 0.16 ) node[] {-};
\draw ( -0.12 , 0 )--( -0.12 , 0.14 ) node[] {-};
\draw ( -0.06 , 0 )--( -0.06 , 0.12 ) node[] {-};
\draw ( 0.0 , 0 )--( 0.0 , 0.1 ) node[] {-};
\draw ( 0.07 , 0 )--( 0.07 , 0.09 ) node[] {-};
\draw ( 0.13 , 0 )--( 0.13 , 0.07 ) node[] {-};
\draw ( 0.2 , 0 )--( 0.2 , 0.06 ) node[] {-};
\draw ( 0.26 , 0 )--( 0.26 , 0.05 ) node[] {-};
\draw ( 0.32 , 0 )--( 0.32 , 0.04 ) node[] {-};
\draw ( 0.39 , 0 )--( 0.39 , 0.04 ) node[] {-};
\draw ( 0.45 , 0 )--( 0.45 , 0.03 ) node[] {-};
\draw ( 0.52 , 0 )--( 0.52 , 0.02 ) node[] {-};
\draw ( 0.58 , 0 )--( 0.58 , 0.02 ) node[] {-};
\draw ( 0.64 , 0 )--( 0.64 , 0.02 ) node[] {-};
\draw ( 0.71 , 0 )--( 0.71 , 0.01 ) node[] {-};
\draw ( 0.77 , 0 )--( 0.77 , 0.01 ) node[] {-};
\draw ( 0.84 , 0 )--( 0.84 , 0.01 ) node[] {-};
\draw ( 0.9 , 0 )--( 0.9 , 0.01 ) node[] {-};
\draw ( 0.97 , 0 )--( 0.97 , 0.01 ) node[] {-};
\draw ( 1.03 , 0 )--( 1.03 , 0.01 ) node[] {-};
\draw ( 1.09 , 0 )--( 1.09 , 0.0 ) node[] {-};
\draw ( 1.16 , 0 )--( 1.16 , 0.0 ) node[] {-};
\draw ( 1.22 , 0 )--( 1.22 , 0.0 ) node[] {-};
\draw ( 1.29 , 0 )--( 1.29 , 0.0 ) node[] {-};
\draw ( 1.35 , 0 )--( 1.35 , 0.0 ) node[] {-};
\draw ( 1.41 , 0 )--( 1.41 , 0.0 ) node[] {-};
\draw ( 1.48 , 0 )--( 1.48 , 0.0 ) node[] {-};
\draw ( 1.54 , 0 )--( 1.54 , 0.0 ) node[] {-};
\draw ( 1.61 , 0 )--( 1.61 , 0.0 ) node[] {-};
\draw ( 1.67 , 0 )--( 1.67 , 0.0 ) node[] {-};
\draw ( 1.73 , 0 )--( 1.73 , 0.0 ) node[] {-};
\draw ( 1.8 , 0 )--( 1.8 , 0.0 ) node[] {-};
\draw ( 1.86 , 0 )--( 1.86 , 0.0 ) node[] {-};
\draw ( 1.93 , 0 )--( 1.93 , 0.0 ) node[] {-};
\draw (0,0) node[below] {\footnotesize 0};
\draw (-1,0) node[below] {\footnotesize -1};
\draw (-2,0) node[below] {\footnotesize -2};
\draw (-3,0) node[below] {\footnotesize -3};
\draw (-4,0) node[below] {\footnotesize -4};
\draw (-5,0) node[below] {\footnotesize -5};
\draw (-6,0) node[below] {\footnotesize -6};
\draw (-7,0) node[below] {\footnotesize -7};
\draw (-8,0) node[below] {\footnotesize -8};
\draw[->] (2.25,0) -- (2.25,3.25);
\draw (2.25,1) -- (2.4,1);
\draw (2.25,2) -- (2.4,2);
\draw (2.25,3) -- (2.4,3);
\draw (2.5,1) node[right] {\footnotesize 6000};
\draw (2.5,2) node[right] {\footnotesize 12000};
\draw (2.5,3) node[right] {\footnotesize 18000};
\end{tikzpicture}
$$
}
\caption{A histogram of  $\sum_{N=1}^{\mathcal{N}} \left(\sum_{n=1}^{N} \frac{B(2,n)}{n}-\zeta(2) N+\frac{1}{2}\, \ln(N)+\mu\right)$
for $\mathcal{N}$ up to $10^6$. There are 250 bins.\label{fig:hist}}
\end{figure}

Another short Python script to generate a histogram of 
$$
\sum_{N=1}^{\mathcal{N}} \left(\sum_{n=1}^{N} \frac{B(2,n)}{n}-\zeta(2) N+\frac{1}{2}\, \ln(N)+\mu\right)\, ,\ \text{ for $\mathcal{N}$ up to $10^6$,}
$$
is as follows:
\begin{verbatim}
import sympy

gamma_dbl =  np.double(sympy.EulerGamma)
zeta_2_dbl = np.double(sympy.zeta(2))
mu = gamma_dbl/2+np.log(24*zeta_2_dbl)/4-zeta_2_dbl/2
a = plt.hist(np.cumsum(err_data+mu),250)
\end{verbatim}
A plot of the histogram is given in Figure \ref{fig:hist}.

\section{Calculations and inferences for the conjectures}

We may use equation (\ref{eq:Cauchy}).
Suppose that we start from Theorem \ref{thm:Cesaro}.
Then we may approximate
$$
\sum_{n=1}^{\infty} \frac{B(2,n)}{n}\, z^n\, =\, (1-z)\, \sum_{N=1}^{\infty}
\sum_{n=1}^{N} \frac{B(2,n)}{n} z^N\, .
$$
Namely, we approximate it by $z \zeta(2)/(1-z)$.
Using this in equation (\ref{eq:Cauchy})
and using the standard saddle point method of Hayman leads to the approximation
$$
	\frac{A(\ell,n,k)}{n!}\, \stackrel{?}{\sim}\, 
	\frac{1}{2\pi n}\, e^{2k-2k\ln(k)+k\ln(n)+k \ln(\zeta(2))-(k^2/(2n))} (1+o(1))\, .
$$
See, for example, Flajolet and Sedgewick \cite{FlajoletSedgewick}.

\begin{remark}
We write $\stackrel{?}{\sim}$ to indicate that the calculation is at least partly heuristic.
(In the above case because we are basing that calculation on the heuristic replacement of $\sum_{n=1}^{N} \frac{B(2,n)}{n}$
by $N\zeta(2)$.)
Here we are trying to justify why we make the conjecture as it is stated: it is based on heuristic
calculations. The numerics is shown as justification.
\end{remark}

Let us specialize to $k \sim \sqrt{\zeta(2)n}$ which is the maximal value.
Then we see that relative to the anticipated value in Corollary \ref{cor:CLT},
we have an extra factor $e^{-\zeta(2)/2}$.
That is one thing to note.
But we also want to improve the anticipated formula to 
the formula in Conjecture \ref{conj:LCLT}.
Firstly, to get the logarithm term from (\ref{eq:varepsilon}), we may guess
$$
\sum_{n=1}^{N} \frac{B(2,n)}{n}\, \stackrel{?}{\sim}\, N \zeta(2) + c \ln(N)\, ,
$$
for some constant $c$.
By an application of the Euler-Maclaurin sum-to-integral comparison one may 
see that the effect this has is to alter the formula to 
$$
	\frac{A(\ell,n,k)}{n!}\, \stackrel{?}{\sim}\, 
	\frac{1}{2\pi n}\, e^{2k-2k\ln(k)+k\ln(n)+k \ln(\zeta(2))}e^{-\zeta(2)/2}
e^{c\left(\ln(n/k)-\gamma\right)} (1+o(1))\, .
$$
Using $k \sim \sqrt{\zeta(2)/n}$ again, we see that we should choose $c=-1/2$
in which case we get
$$
	\frac{A(\ell,n,k)}{n!}\, \stackrel{?}{\sim}\, 
	\frac{1}{2\pi n}\, e^{2k-2k\ln(k)+k\ln(n)+k \ln(\zeta(2))}e^{-\zeta(2)/2}
e^{\ln(k/n)/2} e^{\gamma/2} (1+o(1))\, .
$$
But in order to get the non-logarithm term of (\ref{eq:varepsilon})
we wanted $e^{\ln(24 \zeta(2))/4}$.
So that shows that we should take 
$$
\sum_{n=1}^{N} \frac{B(2,n)}{n}\, \stackrel{?}{\sim}\, 
N \zeta(2) - \frac{1}{2}\, \ln(N) - \mu\, ,
$$
with $\mu$ of the form written in Conjecture \ref{conj:main}.
Note that we do not make the conjecture precisely as it is written above.
Numerically, there are order-1 fluctuations around this quantity.
But we may also write
\begin{equation} 
\label{eq:doublesum}
\sum_{n=1}^{\infty} \frac{B(2,n)}{n}\, z^n\, =\, (1-z)\, \sum_{N=1}^{\infty}
\sum_{n=1}^{N} \frac{B(2,n)}{n}\, z^N\, = (1-z)^2\, \sum_{\mathcal{N}=1}^{\infty} \sum_{N=1}^{\mathcal{N}} \sum_{n=1}^{N} \frac{B(2,n)}{n}\, z^{\mathcal{N}}\, .
\end{equation}
Then we see that the same conclusion should hold -- namely Conjecture \ref{conj:LCLT} for $\ell=2$ --  as long as Conjecture \ref{conj:main} is true as it is written.
In Figure \ref{fig:hist} we plotted the histogram for
$$
\sum_{N=1}^{\mathcal{N}} \left(\sum_{n=1}^{N} \frac{B(2,n)}{n}-\zeta(2) N+\frac{1}{2}\, \ln(N)+\mu\right)\, ,\ \text{ for $\mathcal{N}$ up to $10^6$.}
$$
Let us call the random variable $\mathsf{X}_{\mathcal{N}}$.
But note that the random variable would be multiplied by $(1-z)^2$ in the generating function by (\ref{eq:doublesum}):
\begin{equation*}
\begin{split}
&-\ln(k!) - n\, \ln(r)+
 k\, \ln\left((1-r)^2\, \sum_{\mathcal{N}=1}^{\infty} \sum_{N=1}^{\mathcal{N}} \sum_{n=1}^{N} \frac{B(2,n)}{n}\, r^{\mathcal{N}}\right)\\
&\hspace{2cm} \stackrel{?}{=}\, 2k  -2k\ln(n) + k \ln(n) + k \ln(\zeta(2)) + \frac{1}{2}\, \ln\left(\frac{k}{n}\right) 
+ \frac{1}{4}\, \ln\left(24 \zeta(2)\right)\\
&\hspace{3cm}
+k\, \ln\left(1+\frac{(1-r)^3}{r \zeta(2)}\, \sum_{\mathcal{N}=1}^{\infty}  r^{\mathcal{N}}\, \mathsf{X}_{\mathcal{N}}\right) + o(1)\, .
\end{split}
\end{equation*}
for $r=1-(k/n)$ and $k \sim \sqrt{\zeta(2)n}$ as $n \to \infty$.
Replacing $\mathsf{X}_{\mathcal{N}}$ by an order-1 quantity equal to its supposed limiting theoretical mean -- $\mathbf{E}[\mathsf{X}] = \lim_{\mathfrak{N} \to \infty}
\frac{1}{\mathfrak{N}} \sum_{\mathcal{N}=1}^{\mathfrak{N}}\mathsf{X}_{\mathcal{N}}$ -- 
would still lead only to a quantity $\mathbf{E}[\mathsf{X}]/\zeta(2)$ times $k^3/n^2$ which is $O(1/\sqrt{n})$: in particular $o(1)$.

\appendix
\section{The probability model}
\label{sec:q}
Let us first note that by the $q$-binomial formula
$$
(z;q)_{\ell-1}\, =\, \sum_{j=0}^{\ell-1} q^{j(j-1)/2}\, \frac{(q;q)_{\ell-1}}{(q;q)_{j}(q;q)_{\ell-1-j}}
(-z)^j\, .
$$
For instance, see Section 10.2 of \cite{AndrewsAskeyRoy}.
Therefore,
$$
(q^{k+1}z;q)_{\ell-1}\, =\, \sum_{j=0}^{\ell-1} q^{j(j-1)/2}\, \frac{(q;q)_{\ell-1}}{(q;q)_{j}(q;q)_{\ell-1-j}}
(-z)^jq^jq^{kj}\, .
$$
Therefore,
$$
	\sum_{k=0}^{\infty} q^k (q^{k+1}z;q)_{\ell-1}\,
	=\, \sum_{j=0}^{\ell-1} q^{j(j-1)/2}\, \frac{(q;q)_{\ell-1}}{(q;q)_{j}(q;q)_{\ell-1-j}}
(-z)^jq^j\sum_{k=0}^{\infty} q^{(j+1)k}\, .
$$
So, by the geometric series formula
$$
	\sum_{k=0}^{\infty} q^k (q^{k+1}z;q)_{\ell-1}\,
	=\, \sum_{j=0}^{\ell-1} q^{j(j-1)/2}\, \frac{(q;q)_{\ell-1}}{(q;q)_{j}(q;q)_{\ell-1-j}}
(-z)^jq^j\, \frac{1}{1-q^{j+1}}\, .
$$
But $(q;q)_j (1-q^{j+1}) = (q;q)_{j+1}$.
Shifting the summation index in the identity,
$$
	\sum_{k=0}^{\infty} q^k (q^{k+1}z;q)_{\ell-1}\,
	=\, \sum_{j=0}^{\ell-1} q^{j(j-1)/2}\, \frac{(q;q)_{\ell-1}}{(q;q)_{j+1}(q;q)_{\ell-1-j}}
(-z)^jq^j\, ,
$$
and using the $q$-binomial formula again, gives the first part of Theorem \ref{thm:Cesaro}.

Now, given a number $n \in \N$, let us define $a_1(n),a_2(n),\dots \in \{0,1,\dots\}$ such that
the prime factorization is 
$$
n\, =\, p_1^{a_1(n)} p_2^{a_2(n)} \cdots \, .
$$
According to the probability model, for any finite number of fixed primes $p_r$ for $r=1,\dots,R$,
if we take the uniform measure on $n \in \{1,\dots,N\}$ in the large $N$ limit, the joint distribution of $a_1,\dots,a_R$
converges weakly to independent-product measure and such that $a_r$ is Geometric-$p_r^{-1}$.
Now let us prove the second part of the theorem.

\begin{proofof}{\bf Proof of Theorem \ref{thm:Cesaro}:}
We first note that for any fixed $r$, we have
$$
	\lim_{N \to \infty} \frac{1}{N}\, \sum_{n=1}^{N} \frac{B({\ell},p_r,a_r(n))}{p_r^{(\ell-1)a_r(n)}}\,
=\, -(1-p_r^{-1}) \left(\frac{(1;p_r^{-1})_{\ell}}{(p_r^{-1};p_r^{-1})_\ell} - \frac{1}{(p^{-1};p^{-1})_{\ell}}\right)\, ,
$$
by the $q$-identity we just established.

But we note that $(1;p_r^{-1})_{\ell} = (1-1) \cdot (p_r^{-1};p_{r}^{-1})_{\ell-1}=0$.
So we have
$$
	\lim_{N \to \infty} \frac{1}{N}\, \sum_{n=1}^{N} \frac{B(\ell,p_r,a_r(n))}{p_r^{(\ell-1)a_r(n)}}\,
=\, \frac{1}{(1-p_r^{-2})\cdots (1-p_r^{-\ell})}\, .
$$
Now, taking a product over $r\leq R$ we get, by the results in \cite{ATB}, that
$$
	\lim_{N \to \infty} \frac{1}{N}\, \sum_{n=1}^{N} \prod_{r=1}^{R}\frac{B(\ell,p_r,a_r(n))}{p_r^{(\ell-1)a_r(n)}}\,
=\, \prod_{r=1}^{R}\frac{1}{(1-p_r^{-2})\cdots (1-p_r^{-\ell})}\, .
$$
By the monotone convergence theorem, we have
$$
	\lim_{N \to \infty} \frac{1}{N}\, \sum_{n=1}^{N} \frac{B(\ell,n)}{n^{\ell-1}}\,
=\, \prod_{r=1}^{\infty}\frac{1}{(1-p_r^{-2})\cdots (1-p_r^{-\ell})}\, =\, \zeta(2)\cdots\zeta(\ell)\, ,
$$
by the Euler product formula. (See for example Section 8.1 of \cite{SteinShakarchi}.)
\end{proofof}

\section{Discrete tori with twists}

\label{app:tori}

Suppose $\ell \in \N$ and we have numbers $f_1,\dots,f_{\ell} \in \N$
and the product is $n=f_1\cdots f_{\ell}$.
Given extra $r$-tuples
$$
\boldsymbol{\phi}_r = (\phi_r^{(1)},\dots,\phi_r^{(r)}) \in [f_1]\times \cdots \times [f_r]\, ,
$$
for $r=1,\dots,\ell-1$, we define a graph which we denote
$$
	\mathbf{T}^{(\ell)}_{f_1,\dots,f_{\ell}}(\boldsymbol{\phi}_1,\dots,\boldsymbol{\phi}_{\ell-1})\, .
$$
Firstly, take the vertex set to be $[f_1]\times \cdots \times [f_{\ell}]$.
Secondly, considering this vertex set to be a subset of the vertices of $\Z^{\ell}$
with the usual edge structure, our graph will contain all the edges of the induced
graph of $[f_1]\times \cdots \times [f_{\ell}]$.
But it will also contain edges joining opposite faces.
For the first direction, every pair of vertices
$$
(f_1,i_2,\dots,i_{\ell})\ \text{ and }\
(1,i_2,\dots,i_{\ell})
$$
are joined by an edge.
Let us define a bijection $\pi_1$ of $[f_1]\times \cdots \times [f_{\ell}]$
$$
\pi_1((i_1,i_2,\dots,i_{\ell}))\, =\, \begin{cases} (i_1+1,i_2,\dots,i_{\ell}) & \text{ if $i_1 < f_1$, and}\\
(1,i_2,\dots,i_{\ell}) & \text{ if $i_1=f_1$.}
\end{cases}
$$

But for the second coordinate we make use of $\boldsymbol{\phi}_1 = (\phi_1^{(1)})$
for $\phi_1^{(1)} \in [f_1]$. Namely, let $\rho_1$ be a bijection of $[f_1]$ wherein
$$
\rho(i)\, =\, \begin{cases} i+1 & \text{ if $i< f_1$, and}\\
i+1-f_1 & \text{ if $i\geq f_1$.}
\end{cases}
$$
So this is addition by $1$  modulo $f_1$ defined on $[f_1]$ (instead of on
$\{0,\dots,f_1-1\}$ as is more typical).
Moreover, $\pi_1(i_1,i_2,\dots,i_{\ell}) = (\rho_1(i_1),i_2,\dots,i_{\ell})$.

Now, for the second direction, every pair of vertices
$$
(i_1,f_2,i_3,\dots,i_{\ell})\ \text{ and }\
(\rho_1^{\phi_1-1}(i_1),1,i_3,\dots,i_{\ell})
$$
are joined by an edge.
Now let us define a bijection $\rho_2$ of $[f_1]\times [f_2]$ defined by the mapping
$$
\rho_2((i_1,i_2))\, =\, \begin{cases}
(i_1,i_2+1) & \text{ if $i_2<f_2$, and}\\
(\rho_1^{\phi_1-1}(i_1),1) & \text{ if $i_2=f_2$.}
\end{cases}
$$
And let us define $\pi_2$ the bijection of $[f_1]\times \cdots \times [f_{\ell}]$
such that $\pi_2(i_1,i_2,i_3,\dots,i_{\ell}) = (\rho_2(i_1,i_2),i_3,\dots,i_{\ell})$
(putting in and taking away parentheses wherever that makes sense).

We note that for each $(\phi_1,\phi_2) \in [f_1] \times [f_2]$ we have
$(i_1,i_2) = \rho_2^{\phi_2-1}(\rho_1^{\phi_1-1}(1),1)$.
It is also easy to see that $\pi_1$ and $\pi_2$ commute because $\rho_2$
used the $\rho_1$ for the twist between the two faces.

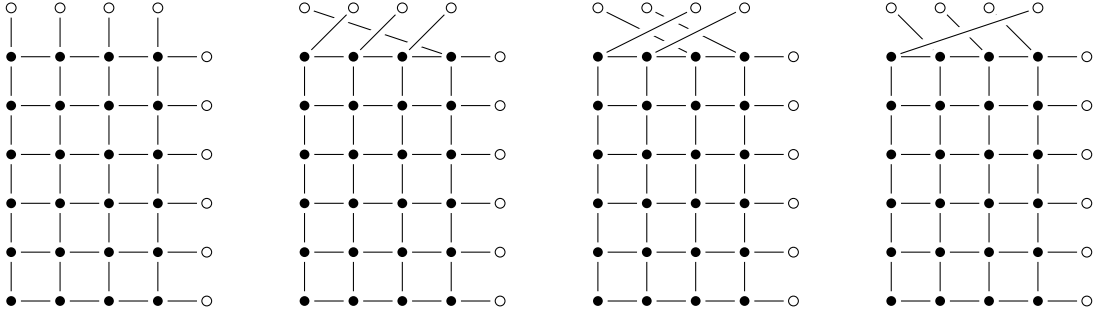
\begin{figure}
\begin{center}
\begin{tikzpicture}[xscale=0.65,yscale=0.65]
\begin{scope}[xshift=-6cm]
\foreach \x in {0,1,2,3}
	{\draw (\x,5) -- +(0,1);}
\foreach \y in {0,1,2,3,4,5}
	{\draw (0,\y) -- (4,\y);
	 \fill[white] (4,\y) circle (2mm);
	 \draw (4,\y) circle (1mm);}
\foreach \x in {0,1,2,3}
	{	\draw (\x,0) -- (\x,5);
		\foreach \y in {0,1,2,3,4,5}
			{   \fill[white] (\x,\y) circle (2mm);
				\fill (\x,\y) circle (1mm);
			}
		\fill[white] (\x,6) circle (2mm);
		\draw (\x,6) circle (1mm);
	}
\end{scope}
\draw (3,5) -- (0,6);
\foreach \x in {0,1,2}
	{\draw[line width=5pt,white] (\x,5) -- +(1,1);
	\draw (\x,5) -- +(1,1);}
\foreach \y in {0,1,2,3,4,5}
	{\draw (0,\y) -- (4,\y);
	 \fill[white] (4,\y) circle (2mm);
	 \draw (4,\y) circle (1mm);}
\foreach \x in {0,1,2,3}
	{	\draw (\x,0) -- (\x,5);
		\foreach \y in {0,1,2,3,4,5}
			{   \fill[white] (\x,\y) circle (2mm);
				\fill (\x,\y) circle (1mm);
			}
		\fill[white] (\x,6) circle (2mm);
		\draw (\x,6) circle (1mm);
	}
\begin{scope}[xshift=6cm]
\draw (2,5) -- (0,6);
\draw (3,5) -- (1,6);
\foreach \x in {0,1}
	{\draw[line width=5pt,white] (\x,5) -- +(2,1);
	  \draw (\x,5) -- +(2,1);}
\foreach \y in {0,1,2,3,4,5}
	{\draw (0,\y) -- (4,\y);
	 \fill[white] (4,\y) circle (2mm);
	 \draw (4,\y) circle (1mm);}
\foreach \x in {0,1,2,3}
	{	\draw (\x,0) -- (\x,5);
		\foreach \y in {0,1,2,3,4,5}
			{   \fill[white] (\x,\y) circle (2mm);
				\fill (\x,\y) circle (1mm);
			}
		\fill[white] (\x,6) circle (2mm);
		\draw (\x,6) circle (1mm);
	}
\end{scope}
\begin{scope}[xshift=12cm]
\draw (1,5) -- (0,6);
\draw (2,5) -- (1,6);
\draw (3,5) -- (2,6);
\draw[line width=5pt,white] (0,5) -- +(3,1);
\draw (0,5) -- +(3,1);
\foreach \y in {0,1,2,3,4,5}
	{\draw (0,\y) -- (4,\y);
	 \fill[white] (4,\y) circle (2mm);
	 \draw (4,\y) circle (1mm);}
\foreach \x in {0,1,2,3}
	{	\draw (\x,0) -- (\x,5);
		\foreach \y in {0,1,2,3,4,5}
			{   \fill[white] (\x,\y) circle (2mm);
				\fill (\x,\y) circle (1mm);
			}
		\fill[white] (\x,6) circle (2mm);
		\draw (\x,6) circle (1mm);
	}
\end{scope}
\end{tikzpicture}
\end{center}
\caption{Four discrete $\ell=2$ dimensional tori with a twist. We have $n=24$, $f_1=4$, $f_2=6$.
Reading left-to-right the choices of $\phi^{(2)}_1 \in [f_1]$ are $1,2,3,4$.
The white circles represent the vertices on the opposite face, repeated so as to 
draw the graphs without extraneous crossings.
\label{fig:torus}}
\end{figure}

\begin{figure}
\begin{center}
\begin{tikzpicture}[xscale=0.9,yscale=0.8]
\draw (2,5) -- (0,6);
\draw (3,5) -- (1,6);
\foreach \x in {0,1}
	{\draw[line width=5pt,white] (\x,5) -- +(2,1);
	\draw (\x,5) -- +(2,1);}
\foreach \y in {0,1,2,3,4,5}
	{\draw (0,\y) -- (4,\y);
	 \fill[white] (4,\y) circle (2mm);
	 \draw (4,\y) circle (1mm);}
\foreach \x in {0,1,2,3}
	{	\draw (\x,0) -- (\x,5);
		\foreach \y in {0,1,2,3,4,5}
			{   \fill[white] (\x,\y) circle (2mm);
				\fill (\x,\y) circle (1mm);
			}
		\fill[white] (\x,6) circle (2mm);
		\draw (\x,6) circle (1mm);
	}
\begin{scope}[xshift=6cm,yshift=-0.75cm,xscale=1.5,yscale=1.5,red]
\draw[line width=4pt,white]  (0,5) -- +(3,-4);
		\draw[->>>>>>>] (0,5) -- +(3,-4);
\draw[line width=4pt,white]  (1,5) -- +(-1,-4);
		\draw[->>>>>>>] (1,5) -- +(-1,-4);
\draw[line width=4pt,white]  (2,5) -- +(-1,-4);
		\draw[->>>>>>>] (2,5) -- +(-1,-4);
\draw[line width=4pt,white]  (3,5) -- +(-1,-4);
		\draw[->>>>>>>] (3,5) -- +(-1,-4);
\draw[blue,line width=4pt,white]  (0,4) -- +(3,-4);
		\draw[blue,->>>>>>>] (0,4) -- +(3,-4);
\draw[blue,line width=4pt,white]  (1,4) -- +(-1,-4);
		\draw[blue,->>>>>>>] (1,4) -- +(-1,-4);
\draw[blue,line width=4pt,white]  (2,4) -- +(-1,-4);
		\draw[blue,->>>>>>>] (2,4) -- +(-1,-4);
\draw[blue,line width=4pt,white]  (3,4) -- +(-1,-4);
		\draw[blue,->>>>>>>] (3,4) -- +(-1,-4);
\foreach \y in {0,1,2,3}
		{\draw[line width=4pt,white]  (3,\y) -- +(-3,2);
		\draw[->>>>>>>] (3,\y) -- +(-3,2);}
\draw[blue,line width=4pt,white]  (3,0) -- +(-3,2);
		\draw[blue,->>>>>>>] (3,0) -- +(-3,2);
\draw[blue,line width=4pt,white]  (3,2) -- +(-3,2);
		\draw[blue,->>>>>>>] (3,2) -- +(-3,2);
\foreach \x in {0,1,2}
	{\foreach \y in {0,1,2,3}
		{\draw[line width=4pt,white]  (\x,\y) -- +(1,2);
		\draw[->>>>>>>] (\x,\y) -- +(1,2);}}
\draw[line width=4pt,white]  (0,0) -- +(1,2);
\draw[blue,->>>>>>>] (0,0) -- +(1,2);
\draw[line width=4pt,white]  (1,2) -- +(1,2);
\draw[blue,->>>>>>>] (1,2) -- +(1,2);
\draw[line width=4pt,white]  (1,0) -- +(1,2);
\draw[blue,->>>>>>>] (1,0) -- +(1,2);
\draw[line width=4pt,white]  (2,2) -- +(1,2);
\draw[blue,->>>>>>>] (2,2) -- +(1,2);
\draw[line width=4pt,white]  (2,0) -- +(1,2);
\draw[blue,->>>>>>>] (2,0) -- +(1,2);
\draw[line width=4pt,white]  (0,2) -- +(1,2);
\draw[blue,->>>>>>>] (0,2) -- +(1,2);
\foreach \x in {0,1,2,3}
	{\foreach \y in {0,1,2,3,4,5}
			{   \fill[white] (\x,\y) circle (1.5mm);
				\fill[black] (\x,\y) circle (1mm);
			}
	}
\end{scope}
\end{tikzpicture}
\end{center}
\caption{An indication of a further twist in case $\ell=3$. Suppose we choose $f_1=4$, $f_2=6$
and $f_3=n/(f_1f_2)$ for some number $n \in f_1 f_2 \N$.
We take $\phi^{(2)}_1 = 4 \in [f_1]$.
Now let us take $(\phi^{(3)}_1,\phi^{(3)}_2) = (2,3) \in [f_1] \times [f_2]$.
Then we give an indication of the bijection $\pi$ on $[f_1]\times[f_2]$ such that to go from $(x,y,f_3)$
to a point on the opposite face we use $(\pi(x,y),1)$.
The permutation is shown in cycle decomposition. There are two cycles for this particular permutation:
one shown in blue, one in red.
\label{fig:torus2}
}
\end{figure}
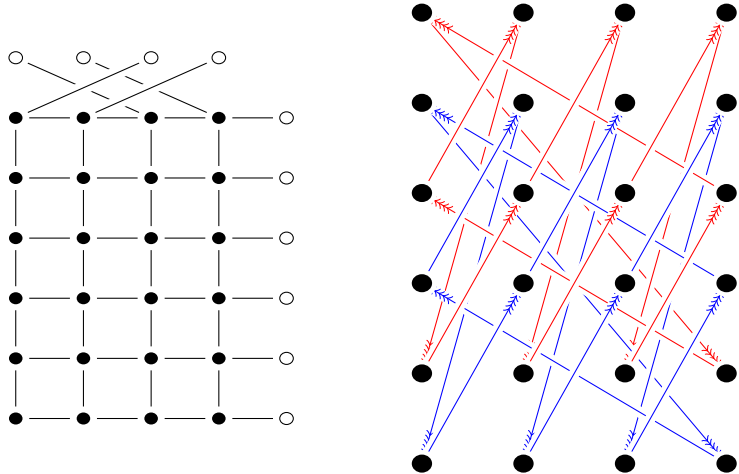

Now we proceed inductively.
Suppose we have bijections $\pi_1,\dots,\pi_{r-1}$ coming from $\rho_1,\rho_2,\dots,\rho_{r-1}$. Given $(\phi_{r-1}^{(1)},\dots,\phi_{r-1}^{(r-1)})$
choose the pair of vertices
$$
(i_1,\dots,i_{r-1},f_r,i_{r+1},\dots,i_{\ell})\ \text{ and }\ 
(j_1^{(r-1)},\dots,j_{r-1}^{(r-1)},1,i_{r+1},\dots,i_{\ell})
$$
to be connected by an edge, where we define
$$
j_1^{(1)} = \rho_1^{\phi_r^{(1)}-1}(i_1)
$$
then
$$
(j_1^{(2)},j_2^{(2)}) = \rho_2^{\phi_r^{(2)}-1}(j_1^{(1)},i_2)
$$
then
$$
(j_1^{(3)},j_2^{(3)},j_3^{(3)}) = \rho_3^{\phi_r^{(3)}-1}(j_1^{(2)},j_2^{(2)},i_3)\, ,
$$
and so on up to $(j_1^{(r-1)},\dots,j_{r-1}^{(r-1)})$.
Note that if $(i_1,\dots,i_{r-1})=(1,\dots,1)$ then we do not cross any 
edges so that $(j_1^{(r-1)},\dots,j_{r-1}^{(r-1)}) = (\phi_1,\dots,\phi_{r-1})$.
In Figure \ref{fig:torus}, we give examples for $\ell=2$.
For $\ell=3$ it is harder to visualize the examples.
But in Figure \ref{fig:torus2} we show an example of the ``twist,''
meaning the permutation we use in going from the $i_3=f_3$ face back to the $i_3=1$ face.

Then we define 
$$
\pi_r(i_1,\dots,i_{\ell})\, =\, 
\begin{cases}
	(i_1,\dots,i_{r-1},i_r+1,i_{r+1},\dots,i_{\ell}) & \text{ if $i_r<f_r$, and}\\
	(j_1^{(r-1)},\dots,j_{r-1}^{(r-1)},1,i_{r+1},\dots,i_{\ell}) & \text{ if $i_r=f_r$.}
\end{cases}
$$
Bryan and Fulman made use of wreath products.

Inductively, $|\langle \pi_1,\dots,\pi_{\ell} \rangle|$ equals $n$: every element may be written as $\pi_1^{\phi_1-1} \cdots \pi_{\ell}^{\phi_{\ell}-1}$
for $(\phi_1,\dots,\phi_{\ell}) \in [f_1]\times \cdots \times [f_{\ell}]$.
Then, by a standard double-counting argument, as in
$$
\text{\url{https://en.wikipedia.org/wiki/Double_counting_(proof_technique)}}
$$
we see directly why $|A(\ell,n,1)|$ satisfies the formula written by Abdesselam.

\section*{Acknowledgments}
I am grateful to Malek Abdesselam.


\baselineskip=12pt
\bibliographystyle{plain}

\end{document}